\documentclass[12pt]{amsart}
\usepackage{amsfonts, amssymb,amsmath, amsthm, amsxtra, latexsym, amscd}
\usepackage[latin1]{inputenc} % for the swedish symbols äåö etc
\def\draft{\centerline{(Draft {\the \day}/{\the\month} \the \year.)}}

\theoremstyle{definition}
\newtheorem{theo+}    {Theorem}      [section]
\newtheorem{prop+}  [theo+]  {Proposition}
\newtheorem{coro+}  [theo+]  {Corollary}
\newtheorem{lemm+}  [theo+]  {Lemma}
\newtheorem{deep+}  [theo+]  {Deep Result}
\newtheorem{fact+}  [theo+]  {Fact}
\theoremstyle{definition}
\newtheorem{exam+}  [theo+]  {Example}
\newtheorem{rema+}  [theo+]  {Remark}
\newtheorem{defi+}  [theo+]  {Definition}
\newtheorem{xca+}[theo+]{Exercise}

\numberwithin{equation}{section}

\def\draft{\centerline{(Draft {\the \day}/{\the\month} \the \year.)}}
\bigskip

\def\refn#1.#2{\expandafter\def\csname#1\endcsname{[#2]}}
\def\refnr#1.{\csname#1\endcsname}

\def\Del{\Delta}

\def\a{\alpha}
\def\b{\beta}

\def\Claminv2{|C(\Lambda)|^{-2}}

\def\Ga{\Gamma}
\def\varepsi{\varepsilon}
\def\lam{\lambda}

\def\ome{\omega}
\def\Ome{\Omega}

\def\de{d\varepsilon}

\def\Aa2D{A^{\a,2}(D)}
\def\bAa2D{\overline{A^{\a,2}(D)}}
\def\Ab2D{A^{\beta,2}(D)}
\def\bAb2D{\overline{A^{\beta,2}(D)}}

\def\Norm#1_#2{\Vert#1\Vert_{#2}}

\def\phipl12{\phi_{p_{l_1}, p_{l_2}}}
\def\phip01{\phi_{p_{0}, p_{0}}}

\def\a{\alpha}
\def\b{\beta}

\def\Claminv2{|C(\Lambda)|^{-2}}

\def\Del{\Delta}

\def\varepsi{\varepsilon}

\def\Ga{\Gamma}
\def\lam{\lambda}

\def\det{\operatorname{det}}
\def\diag{\operatorname{diag}}

\def\tr{\operatorname{tr}}

\def\de{d\varepsilon}

\def\Aa2D{A^{\a,2}(D)}
\def\bAa2D{\overline{A^{\a,2}(D)}}
\def\Ab2D{A^{\beta,2}(D)}
\def\bAb2D{\overline{A^{\beta,2}(D)}}

\def\phipl12{\phi_{p_{l_1}, p_{l_2}}}
\def\phip01{\phi_{p_{0}, p_{0}}}

\def\bc{\mathbb C}
\def\br{\mathbb R}

\def\alg/{algebra} %Ê\alg/ic 
\def\Alg/{Algebra} %Ê\Alg/ic 
\def\alt/{alternative} % \alt/ly
\def\anal/{analytic}
\def\analfunc/{\anal/\ \func/}
\def\Ans/{\it Answer. \normal}
\def\ass/{associative}
\def\nass/{non-\ass/}
\def\autom/{automorphism}
\def\homom/{homomorphism}
\def\isom/{isomorphism}
\def\bdd/{bounded}
\def\Bdd/{Bounded}
\def\bddsymdom/{bounded \sym/ \dom/}
\def\Cartdom/{Cartan \dom/}
\def\bdry/{boundary}
\def\bsd/{\bdd/ \symdom/}
\def\bv/{boundary value}
\def\cf/{{\it cf}\.}
\def\Cf/{{\it Cf}\.}
\def\charr/{character}
\def\coeff/{coefficient}
\def\comm/{commutative}
\def\cpct/{compact}
\def\compl/{complex}
\def\comp/{complex}
\def\Comp/{Complex}
\def\conf/{conformal}
\def\conj/{conjugate}
\def\conn/{connect}
\def\cont/{continuous}
\def\conv/{converge} % \conv/nce \conv/nt
\def\convc/{convergence}
\def\convt/{convergent}
\def\convx/{convex}
\def\coord/{coordinate}
\def\lcoord/{local coordinate}
\def\Corr/{Corresponding}
\def\corr/{corresponding}
\def\corrd/{correspond}
\def\cov/{covariant}
\def\decomp/{decomposition}
\def\deco/{decompose}
\def\diff/{different} % \diff/iable \diff/ial
\def\Diff/{Different} % \Diff/able \Diff/ial
\def\dimn/{dimension} % \dimen/al
\def\distr/{distribution} % \distr/al
\def\div/{diverge} % \div/nt 
\def\dom/{domain}
\def\eg/{\hbox{\it e.g}\.}
\def\eigenf/{eigen\-\func/}
\def\eigensp/{eigen\-space}
\def\eigenv/{eigen\-value}
\def\eq/{equation}
\def\equa/{equation}
\def\de/{\diff/ial \equa/}
\def\do/{\diff/ial operator}
\def\ode/{ordinary \de/}
\def\pde/{partial \de/}
\def\pdo/{partial \diff/ial operator}
\def\psdo/{pseudo \diff/ial operator}
\def\fin/{finite}
\def\Ex/{\it Example.\ \normal}
\def\Exnr#1/{\it Example #1.\ \normal}
\def\foll/{follow}
\def\follg/{following}
\def\Follg/{Following}
\def\func/{function}
\def\Func/{Function}
\def\Fonc/{Fonc\-tion}
\def\fonc/{fonc\-tion}
\def\Funk/{Funk\-tion}
\def\funk/{Funk\-tion}
\def\gen/{general}
\def\har/{harmonic}
\def\Hint/{\it Hint. \normal}
\def\hist/{historic}
\def\histcl/{historical}
\def\hol/{holo\-morphic}
\def\homog/{ho\-mo\-ge\-ne\-ous}
\def\hyp/{hyper\-bolic}
\def\hyperg/{hyper\-geometric}
\def\ie/{\hbox{\it i.e.}}
\def\iff/{if and only if}
\def\ineq/{inequality}
\def\infra/{{\it inf\-ra}}
\def\ultra/{{\it ult\-ra}}
\def\Inpart/{In particular}
\def\inpart/{in particular}
\def\instof/{instead of}
\def\interps/{interpolation space}
\def\interp/{interpolation}
\def\Interp/{Interpolation}
\def\interpr/{Interpretation}
\def\Intr/{Introduction}
\def\intv/{interval}
\def\inv/{invariant}
\def\invc/{invariance}
\def\Iowords/{In other words}
\def\iowords/{in other words}
\def\ipr/{inner product}
\def\irred/{irreducible}
\def\lb/{line bundle}
\def\lin/{linear}
\def\lhs/{left hand side}
\def\rhs/{right hand side}
\def\loc/{local}
\def\math/{mathematic} 
\def\mathcn/{\math/ian}
\def\manif/{manifold}
\def\meas/{measure}
\def\measl/{measurable}
\def\mero/{mero\-morphic}
\def\mon/{monomial}
\def\monog/{monogenic}
\def\mult/{multiple}
\def\multy/{multiply}
\def\multn/{multiplication}
\def\nas/{necessary and sufficient}
\def\nbd/{neighborhood}
\def\neg/{negative}
\def\nondeg/{nondegenerate}
\def\Oohand/{On the other hand}
\def\oohand/{on the other hand}
\def\Oonhand/{On the one hand}
\def\oonhand/{on the one hand}
\def\oper/{operator}
\def\orth/{ortho\-gonal}
\def\orthon/{ortho\-normal}
\def\otoh/{on the other hand}
\def\quat/{quaternion}
\def\pp/{\hbox{a. e.}}
\def\psh/{plurisubharmonic}
\def\pol/{polynomial}
\def\pot/{potential}
\def\pos/{positive}
\def\princ/{principle}
\def\prob/{probability}
\def\proj/{projective}
\def\projn/{projection}
\def\Proof/{\it Proof:\normal}
\def\Rem/{\it Remark\normal}
\def\Remnr#1/{\it Remark\ \normal #1. }
\def\rep/{representation}
\def\meta/{metaplectic representation}
\def\repr/{reproducing}
\def\reprker/{reproducing kernel}
\def\resp/{respective} % \resp/ly
\def\resply/{respectively}
\def\restr/{restriction}
\def\sa/{self-adjoint}
\def\st/{such that}

\def\sol/{solution}
\def\ru/{space}
\def\sph/{spherical}
\def\ssp/{sub\ru/}
\def\sym/{symmetric}
\def\Sym/{Symmetric}
\def\symb/{symbol}
\def\symbc/{symbolic}
\def\symdom/{\sym/ domain}
\def\symp/{symplectic}
\def\Theor#1/{\fet Theorem #1.\ \normal}
\def\Lem#1/{\fet Lemma #1.\ \normal}
\def\Lemma/{\fet Lemma.\ \normal}
\def\topl/{topology}
\def\topll/{topological}
\def\transf/{transform}
\def\transl/{translation}
\def\transfn/{transformation}
\def\transv/{transvectant}
\def\trig/{trigonometric}
\def\tril/{trilinear}
\def\trilf/{trilinear form}
\def\uhp/{upper halfplane}
\def\uhs/{upper halfspace}
\def\vb/{vector bundle}
\def\vf/{vector field}
\def\vsp/{vector space}
\def\wrt/{with respect to}
\def\Wlog/{Without loss of generality}
\def\a{\alpha}

\def\lam{\lambda}

\def\Ab/{Abel}
\def\Ban/{Banach}
\def\Bansp/{\Ban/ space}
\def\Belt/{Bel\-tra\-mi}
\def\Berg/{Berg\-man}
\def\Bern/{Ber\-nou\-lli}
\def\Berz/{Berezin}
\def\Bess/{Bessel}
\def\Cart/{Car\-tan}
\def\Cay/{Cay\-ley}
\def\CG/{Clebsch-Gordan}
\def\Cl/{Clifford}
\def\CR/{Cauchy-Rie\-mann}
\def\Dir/{Dirichlet}
\def\Eucl/{Euclide}
\def\F/{Fourier}
\def\Hank/{Hankel}
\def\Hankf/{\Hank/ form}
\def\Herm/{Hermite}
\def\Hilb/{Hilbert}
\def\Hilbs/{Hilbert space}
\def\Hilbsp/{Hilbert space}
\def\HS/{Hilbert-Schmidt}
\def\Lag/{La\-grange}
\def\Lap/{La\-place}
\def\LapBelt/{\Lap/-\Belt/}
\def\Leb/{Lebesgue}
\def\Marc/{Mar\-cin\-kie\-wicz}
\def\Moeb/{Moebius}
\def\Moebt/{Moebius transformation}
\def\Moebtransfn/{Moebius transformation}
\def\Pla/{Plan\-che\-rel}
\def\Poin/{Poin\-car\'e}
\def\Riem/{Rie\-mann}
\def\Riemn/{\Riem/ian}
\def\psRiemn/{pseudo-\Riem/ian}
\def\Riems/{Rie\-mann surface}
\def\Schroe/{Schr\"odinger}
\def\Weier/{Weier\-strass}
\def\anal/{analytic}
\def\bsd/{bounded symmetric domain  }
\def\bdd/{bounded}
\def\calc/{calculation}\def\conj{conjugate}
\def\calci/{calculating}\def\eg{e.g.}
\def\conj/{conjugate}
\def\deco/{decomposition}
\def\eg/{e.g.}
\def\fct/{function}
\def\gp/{group}
\def\hw/{highest weight}
\def\hwv/{highest weight vector}
\def\hwvs/{highest weight vectors}
\def\lw/{lowest weight}
\def\lwv/{lowest weight vector}
\def\lwvs/{lowest weight vectors}
\def\hds/{holomorphic discrete series}
\def\iff/{if and only if}
\def\inv/{invariant}
\def\irrde/{irreducible decomposition}
\def\meas/{measure}
\def\transf/{transform}
\def\rep/{representation}
\def\resp/{respectively}
\def\inters/{intertwines}
\def\interg/{intertwining}
\def\meta/{metaplectic representation}
\def\qu/{quaternion}
\def\rep/{representation}
\def\symdom/{ symmetric domain}
\def\st/{such that}
\def\shd/{subhead}
\def\transf/{transform}
\def\wrt/{with respect to}

\def\Norm#1#2#3{\Vert#1\Vert^{#3}_{{#2}¥}}

\def\tr{\operatorname{tr}}

\begin{document}
%\draft
\title[Radon transform]
{Radon transform on real, complex
and quaternionic Grassmannians
}
\author{Genkai Zhang}
\address{Department of Mathematics, Chalmers University of Technology and
G\"oteborg University
 , G\"oteborg, Sweden}
\email{genkai@math.chalmers.se}

\thanks{Research  supported by the Swedish
Science Council (VR)}
%\subjclass{????}

\def\bbK{\mathbb K}
\def\bKk{\mathbb K^k}
\def\bKn{\mathbb K^n}
\def\Cos^2{\text{Cos}^2}

\begin{abstract}
Let $G_{n,k}(\bbK)$ be the
 Grassmannian manifold
of $k$-dimensional $\bbK$-subspaces in $\bbK^n$
where $\bbK=\mathbb R, \mathbb C, \mathbb H$
is  the  field of real, complex or quaternionic numbers.
 For $1\le k\le k^\prime \le n-1$ we
 define the Radon transform $(\mathcal R f)(\eta)$, $\eta \in 
G_{n,k^\prime}(\bbK)$, 
for functions $f(\xi)$ on $G_{n,k}(\bbK)$  as an integration over
all $\xi \subset \eta$.
When  $k+k^\prime \le n$ we give an inversion
formula in terms of the G\aa{}rding-Gindikin
fractional integration
and the Cayley type differential operator on
the symmetric cone of positive $k\times k$ matrices
over $\bbK$.
This generalizes the
recent results of Grinberg-Rubin 
for real Grassmannians. 
\end{abstract}

\maketitle



%\curraddr{}


%\author{Author Two}
%\address{}
%\email{}

%    General info

%\date{}

%\dedicatory{}
\def\mP{\mathcal P}
\def\mF{\mathcal F}
\def\SdN{\mathcal S^N}
\def\eSdN{\mathcal S_N^\prime}
\def\Gc{G_{\mathbb C}}
\def\rqs{G\backslash K_{\mathbb C}}

\baselineskip 1.25pc

\section{Introduction}

The Radon transform on rank one symmetric spaces has been studied
extensively and is related to many areas in analysis and geometry;
see
 \cite{Helgason-radonbk}, \cite{Hel2}
for a systematic treatment and e.g. \cite{Grinberg-grs},
 \cite{Kakehi-grs},
 \cite{Strichartz-grs},
\cite{Rubin-adv-math},
\cite{Helgason-rmk-Rubin}
and references therein
for some recent development. For  higher rank symmetric spaces
the theory  is far from been complete
and there has been comparably  less progress. In a remarkable paper
\cite{Grinberg-Rubin} Grinberg and Rubin find an inversion formula
for the Radon transform from functions
on the real Grassmannian $G_{n,k}(\mathbb R)$ of $k$-dimensional
subspaces in $\mathbb R^n$ to functions on $G_{n,k^\prime}(\mathbb R)$,
with $1\le k\le k^\prime\le n-1$, $k+k^\prime \le n$, by using
the Gårding-Gindikin fractional integration on
the space of real symmetric $k\times k$-matrices. It is natural
to ask if the corresponding results hold
for the Grassmannian manifolds over the complex and quaternionic
numbers, and for corresponding  non-compact symmetric spaces
of matrix balls. In the present paper we will answer the question
and prove the  results for complex and quaternionic
Grassmannians. 
We proceed
with a brief summary of our results
 and the technical
tools  needed to prove them.

Let $\mathbb K=\br, \bc, \mathbb H$ be the field of real, complex or quaternionic
numbers with real dimension
$d=1, 2, 4$ and with the standard involution $x\to \bar x$. Let
$G_{n, k}(\mathbb K)$  be the Grassmannian manifold 
of $k$-dimensional subspaces over $\mathbb K$ in $\mathbb K^n$.
For $1\le k\le k^\prime \le n-1$ the Radon
transform  $\mathcal R: C^{\infty}(G_{n, k}(\mathbb K))
\to  C^{\infty}(G_{n, k^\prime}(\mathbb K))$ is defined by
\begin{equation}\label{radon-def-1}
\phi(\eta)=(\mathcal Rf)(\eta)=\int_{\xi\subset \eta}f(\xi)d_\eta \xi, \qquad \eta\in 
G_{n, k^\prime}(\mathbb K),
\end{equation}
where $d_\eta \xi$ is a certain probability
measure on the set $\{\xi\in G_{n, k}(\bbK): \xi\subset \eta\}$ 
invariant with
respect to the group of unitary transformations
of $\eta$; it can be defined using  a
group-theoretic formulation, see (\ref{def-radon}).

When $k=1$ the Grassmannian $G_{n, k}(\mathbb K)$ is the
projective space of all lines $\xi$ in $\mathbb K^n$.
The Radon transform $f(\xi)\to \phi(\eta)$
is inverted by using
the integration of $\phi(\eta)$
 over all subspaces $\eta$ that are at  a fixed
angle $\theta$ with the given $\xi$, and by
a fractional integration \cite{Hel2}
with respect to the variable $\cos^2 \theta$.
Motivated
by that result Grinberg and Rubin
\cite{Grinberg-Rubin} define also
a cosine of an ``angle'' for a $k$-plane $\xi$
and $k^\prime$-plane in $\mathbb R^n$,
defined as a $k\times k$ semi-positive matrix, defined
up unitary $U(k, \bbK)$-equivalence. To be
 more precise we denote
$M_{n,  k}:=M_{n, k}(\bbK)$ the  space 
of $n\times k$-matrices over $\mathbb K$, 
and $S_{n, k}=S_{n, k}(\bbK)$ the  Stiefel manifold 
of all orthogonal $k$-frames in $\bbK^n$, identified
also with the set of all isometric $\bbK$-linear transformations $x\in M_{n, k}:\bbK^k\to \bbK^n$.  The Grassmannian $G_{n, k}$
will be viewed as a quotient space of
$S_{n, k}$ via the mapping $x\in S_{n, k}
\to \xi=\{x\}=x\mathbb K^k\in G_{n, k}$, and we will identify
a function $f(\xi)$ on $G_{n, k}$  as a right $U(k, \bbK)$-invariant
function $f(x)$ on $S_{n, k}$.
We define now, following \cite{Grinberg-Rubin},
the cosine of the ``angle'' $(\eta, x)$
 between the element $x\in S_{n, k}$
and $\eta\in G_{n, k^\prime}(\bbK)$ to be
a semi-positive $k\times k$-matrix,
$$
\Cos^2(\eta, x)=x^\ast P_{\eta}x
$$
where $P_{\eta}$ is the orthogonal projection $\bbK^n \to \eta$.
(One can define an angle ``$(\eta, x)$'' as a $k\times k$-self adjoint matrix
using the spectral calculus
of $t\to \cos^2 t $, however we will only
need $\Cos^2(\eta, x)$.) The matrix
$\Cos^2(\eta, \xi)$ of the ``angle''
between $\eta\in G_{n,k^\prime}$ and $\xi\in G_{n, k}$
can be defined up to the $U(k, \bbK)$ equivalence as
$\Cos^2(\eta, x)$, with $\xi=\{x\}$.

The generalization to higher rank
spaces of the  integration of $\phi(\eta)$
is
$$
(\mathcal T_{r^{\frac 12}}\phi)(\xi):= \int_{U(k, \bbK)} du
\int_{ \Cos^2(\eta, x)=uru^\ast}\phi(\eta) d_{\xi} \eta,
\qquad \xi=\{x\},
$$
where $d_x \eta$ is a group-invariant 
measure on the set of integration.

% Introduce Radon with symbols and state main results
%with m=\frac{a}2+\varepsi

Our  main result in this paper is Theorem 4.6, it
gives an inversion formula expressing $f(\xi)$ in terms
the 
G\aa{}rding-Gindikin fractional integration
$I^\varepsi$ and the Cayley type differential operator
acting on the function
$(\mathcal T_{r^{\frac 12}}\phi)(\xi)$ with respect to 
 the $r$-variable.
This generalizes the result
of Grinberg-Rubin for real Grassmannians \cite{Grinberg-Rubin}.
The principal  technical tools in proving it
are Propositions \ref{bi-stiefel} and
\ref{int-stief},  giving integral formulas 
on matrix spaces and on Stiefel manifolds.
Proposition \ref{bi-stiefel} (ii)
is proved 
for the real Grassmannian 
in \cite{Grinberg-Rubin} by using rather complicated
formulas involving Bessel functions
proved by Herz \cite{Herz}. 
Our idea is instead to find first an integral formula
 on the matrix space, namely Proposition \ref{bi-stiefel} (i), 
and (ii) will be then a direct consequence. This then
answers a question in \cite{Grinberg-Rubin}.
A similar integral formula as our Proposition 3.5 is proved 
in \cite[Lemma 2.5]{Grinberg-Rubin}  using again some interesting
and however tricky computations. Here we use
again integration on the matrix space. We believe
that our proofs are easier, both technically
and conceptually. Nevertheless it will be clear that we are very much inspired
by  that paper.

I would like to thank Prof. Boris Rubin for 
some helpful correspondence
and  for bringing \cite{OR} to my attention.  Thanks are also
due to the referees for the careful reading and comments
on an earlier version of this paper.

\section{Symmetric cones and G\aa{}rding-Gindikin
Fractional Integration}

In this section we fix notations and we recall some known results
on the Gindikin Gamma function on symmetric 
cones; see \cite{FK-book} for a systematic treatment.

Let $\mathbb K=\br, \bc, \mathbb H$ be  the
field of real, complex or quaternionic numbers
with real dimension 
$$
d:=\dim_{\mathbb R} \bbK=1, 2, 4
$$
and let $x\to \bar x$ be the standard conjugation.
Let $M_{n,k}:=M_{n,k}(\mathbb K)$ be the Euclidean space
of all $n\times k$-matrices with entries $x$ in $\mathbb K$.
 For $x\in M_{n,k}$ we denote
 $x^\ast:= \bar x^T$, where $T$ stands for the transpose.
The space $M_{n,k}$ is then a Euclidean space
 with the metric $\tr (xx^\ast)$ and the
corresponding Lebesgue measure $dx$. Throughout
this paper we will fix $1\le k< n$ if nothing
else if specified. The vector space  $\bbK^k$
will be identified as column vectors $v$ with multiplication
by $c\in \bbK$ from right, $v\mapsto vc$,
and all
matrices $x\in M_{n,k}$ will be identified as $\bbK$-linear transformations from $\bbK^k$ to $\bbK^n$ (i.e. $x(vc)=x(v)c$, $v\in \bbK^k$, $c\in\bbK$)
unless something else is specified.

The $\mathbb R$-subspace 
$$\mathcal A=\mathcal A_k:=\{a\in M_{k, k}; a^\ast =a\}
$$
 of the  self-adjoint $M_{k,k}$-matrices forms
a formally real Jordan algebra with the Jordan
product $a\circ b=\frac 12(ab+ba)$ 
and the identity $I:=I_k$. Its (real) dimension
$N$ is given by
$$
\quad N:=N_k=\dim_{\mathbb R} \mathcal A=k+\frac{d}{2} k(k-1)
$$
Let $\Ome=\Ome_k$ be the cone in $\mathcal A$ of positive definite
$k\times k$ $\mathbb K$-matrices.
Let $\Delta$
be the determinant function on $\mathcal A$.
It is a polynomial of degree $k$; for real
and complex matrices it is the usual determinant,
for quaternionic matrices it can be defined
using the Pfaffian by identifying $\mathcal A$ with
a subspace of complex $2k\times 2k$ skew symmetric
complex matrices.  Let $GL_k(\mathbb K)$ be the
group of invertible $k\times k$ $\mathbb K$-matrices. It
acts on $\Ome=\Ome_k$ by $g: r\to g rg^\ast$, and
with this action the cone $\Ome$ becomes
a Riemannian symmetric space, $\Ome=GL_k(\mathbb K)/U(k, \mathbb K)$,
where $U(k, \mathbb K)$ is the orthogonal group
over $\bbK$. The  $GL_k(\mathbb K)$-invariant
measure on $\Ome$ is
\begin{equation}
  \label{eq:inv-m}
d\iota(a):=\Delta(a)^{-N/k}da. 
\end{equation}

The Gindikin Gamma integral is defined and given by
  \begin{equation}
    \label{eq:gamm}
\int_{\Ome} e^{-\tr(s t^\ast)}\Delta(s)^{\lam}
d\iota(s)
=\Gamma_{\Ome}(\lam) \Delta(t)^{-\lam},  \qquad t\in \Ome,  \Re\lam> \frac Nk -1  
  \end{equation}
where
$$
\Gamma_{\Ome}(\lam)=(2\pi)^{(N-k)/2}\prod_{j=1}^k \Gamma(\lam-\frac d2(j-1))
$$
is the Gindikin Gamma function. 
(Note that a different normalization of the measure $ds$ is used
in \cite{Grinberg-Rubin}.) The corresponding
Beta-integral is
\begin{equation}
  \label{eq:beta}
\int_0^s \Del(t)^{\lam-N/k}
\Del(s-t)^{\mu-N/k}dt
=B_\Ome(\lam, \mu)\Del(s)^{\lam +\mu -N/k}, 
\qquad
\Re\lam, \Re\mu> \frac Nk -1  
\end{equation}
with
$$
B_\Ome(\lam, \mu)=\frac{\Ga_{\Ome}(\lam)\Ga_{\Ome}(\mu)}{\Ga_\Ome(\lam +\mu)}.
$$

Let $(0, I)$ be the unit open interval $(0, I)=\{s\in \Omega; s<I\}$.
The G\aa{}rding-Gindikin fractional integral for
a function $f$ on $(0, I)$ is
defined by
$$
(I^{\lam}f)(s)=\frac 1{\Gamma_{\Ome}(\lam)}\int_{0}^s \Delta(s-t)^{\lam-N/k}f(t)dt, 
\quad \Re\lam> \frac Nk , 
$$
where  the integration is understood 
as over the
set $t\in \Ome$, $t< s$. 
For our purpose later (see Lemma 4.3) we
consider the space $L^1((0, I), \Del(I-t)^{\nu-\frac Nk}
 dt)$
with $\alpha\in \mathbb R$. The integral
$(I^{\lam}f)(s)$ is then well-defined for $f\in 
L^1((0, I), \Del(I-t)^{\nu-\frac Nk} dt)$.

Using the previous
formula it is easy to prove the
following semigroup property
$$
I^\lam (I^\mu f)=I^{\lam +\mu}f, \qquad
\Re\lam, \Re\mu> \frac Nk 
$$
for  $f\in 
L^1((0, I), \Del(I-t)^{\nu-\frac Nk}dt)$.

To define the  operator $I^\lam$ 
for smaller $\Re\lam$ we need
the differential
operator $\Delta(\partial)$, the so-called Cayley type differential operator, defined uniquely on $\mathcal A$ 
by requiring that
$$
\Delta(\partial_x) e^{\tr (xy^\ast)}
=\Delta(y)e^{\tr(xy^\ast)}.
$$
It particular it follows from the equality
   (\ref{eq:gamm}) that
$$
\Delta(\partial)\Del(s)^\lam =\prod_{j=1}^k(\lam +\frac d2(j-1))
\Del(s)^{\lam-1}
$$
for any $\lam \in \mathbb C$, which sometimes is referred as
Cayley-Capelli type identity.

We state and prove some elementary results on the integral operator $I^\lam$.
They may have been proved in more general form
in the literature. 

\begin{lemm+} 
\begin{itemize}
\item [(i)]
 Suppose   $\Re\lam>\frac Nk-1$ and $\mu> \frac Nk-1 $. 
The operator $I^\lam$ defines
a bounded operator from 
$L^1((0, I),
\Del(I-t)^{\mu +\Re\lam -\frac Nk}dt)$
into 
$L^1((0, I),
\Del(I-t)^{\mu -\frac Nk}dt).
$
\item [(ii)] Suppose  $\nu > \frac Nk-1 $ and
$f\in 
L^1((0, I),\Del(I-t)^{\nu -\frac Nk}dt)$.
$I^\lam f$  has
an analytic continuation in $\lam\in \mathbb C$ as
a distribution on the space $C_0^\infty(0, I)$
of smooth functions with compact support in $(0, I)$.
\item [(iii)] Suppose 
$\nu > \frac Nk-1 $ and $f\in 
L^1((0, I),\Del(I-t)^{\nu -\frac Nk}dt)$.
Let $m>\frac Nk$ be an integer. 
Then 
$$
\Del(\partial)^m I^m f=f,$$ 
in the sense of
distributions.
\end{itemize}
\end{lemm+}

\begin{proof} (i). Let  $\lam$ and $ \mu$ be as in (i). We
estimate the norm of $I^\lam f$ in
$L^1((0, I),
\Del(I-t)^{\mu -\frac Nk}dt)
$. It is, apart from the constant $\Gamma_\Ome (\lam)^{-1}$,
\begin{equation*}
  \begin{split}
  &\quad \int_0^I |\int_0^s \Delta(s-t)^{\lam-\frac Nk} f(t)dt|
\Delta(I-s)^{\mu -\frac Nk}ds\\
&\le 
\int_0^I \int_0^s \Delta(s-t)^{\Re \lam-\frac Nk} |f(t)|
\Delta(I-s)^{\mu -\frac Nk}dt ds\\
&=\int_0^I \left(\int_t^I \Delta(s-t)^{\Re \lam-\frac Nk} 
\Delta(I-s)^{\mu -\frac Nk}ds\right)|f(t)|dt.
  \end{split}
\end{equation*}
We compute the inner integral. Performing the change
of variables $s=t+P((I-t)^{\frac 12})(v), v\in (0, I)$, where
$x\to P(x)$ is the quadratic representation of 
the Jordan algebra $\mathcal A$ (see \cite{FK-book}).
We have then $ds=\Del(I-t)^{\frac Nk}dv$,
$\Del(I-s)=\Del(I-t)\Del(I-v)$,
$\Del(s-t)=\Del(I-t)\Del(v)$.
That integral is
$$
\Del(I-t)^{\mu+\Re\lam -\frac Nk} 
\int_0^I \Delta(v)^{\Re \lam-\frac Nk} 
\Delta(I-v)^{\mu -\frac Nk}dv
=B(\mu, \Re\lam)\Del(I-t)^{\mu+\Re\lam -\frac Nk}.
$$
The result follows
by substituting this into the previous estimate.

(ii). To define $I^\lam$ for general $\lam$ we note that for
 $\Re\lam>\frac Nk $ and a test function
$\phi\in C_0^\infty(0, I)$, we have,
viewing the operator $\Del(\partial)$ as acting on distributions
\begin{equation*}
  \begin{split}
(\Del(\partial)I^{\lam}f, \phi)
&=(-1)^k (I^{\lam}f, \Del(\partial)\phi)
=(-1)^k\frac 1{\Ga_\Ome(\lam)}\int_0^I (I^\lam f)(s)\Del(\partial)\phi(s) ds
\\
&=(-1)^k\frac 1{\Ga_\Ome(\lam)}\int_0^I 
\left(\int_0^s   \Del(s-t)^{\lam-\frac Nk} f(t)dt\right)
\Del(\partial)\phi(s) ds\\
&=\int_0^I \left((-1)^k\frac 1{\Ga_\Ome(\lam)}
\int_t^I  \Del(s-t)^{\lam-\frac Nk} 
\Del(\partial)\phi(s) ds\right)
 f(t)dt    
  \end{split}
\end{equation*}

To treat the inner integral (for fixed t) we 
change variables $s=t+u, u\in \mathcal A$ and 
denote $\psi(u)=\phi(s)=\phi(t+u)$. 
The function $\psi$ so defined
is then a smooth function on $\mathcal A$ of compact support in the
interval $(-t,  I-t)$,
in particular it is in the
Schwartz space $S(\mathcal A)$. Moreover,
and $\Del(\partial_s)\phi(s)=\Del(\partial_u)\psi(u)$.
The inner integral can be written as
$$
(-1)^k
\frac 1{\Ga_\Ome(\lam)}
\int_0^{I-t}\Del(u)^{\lam-\frac Nk}\Del(\partial)\psi(u) du
=(-1)^k\frac 1{\Ga_\Ome(\lam)}
\int_\Ome\Del(u)^{\lam-\frac Nk}\Del(\partial)\psi(u) du.
$$
This is the Riesz integral as tempered distribution
treated in \cite{FK-book}, VII.
It follows from Theorem VII.2.2(ii) there that
the integral is
$$
\frac 1{\Ga_\Ome(\lam -1)}\int_\Ome\Del(u)^{\lam-1
-\frac Nk}\psi(u) du.
$$
This implies that
$$
(\Del(\partial)I^{\lam}f, \phi)
=(I^{\lam-1}f, \phi).
$$
In this way we can define $I^{\lam}f$ for $\Re \lam  >\frac Nk -2$,
and successively for all $\lam$. Thus $I^\lam$
has analytic continuation for all $\lam \in \mathbb C$. 

(iii) By the same computation we have
$$
(\Del(\partial)^mI^{m}f, \phi)
=\int_0^I\left((-1)^{km}\frac 1{\Ga_\Ome(m)}
\int_\Ome\Del(u)^{m-\frac Nk}\Del(\partial)^m\psi(u) du
\right) f(t)dt
$$ 
and it follows again from 
the first and the third equality in 
Theorem VII.2.2  loc. cit.
that the inner integral is the delta distribution
on $\psi$,
$$(-1)^{km}\frac 1{\Ga_\Ome(m)}
\int_\Ome\Del(u)^{m-\frac Nk}\Del(\partial)^m\psi(u) du
=\delta (\psi) =\psi(0)=\phi(t),
$$
proving (iii).
\end{proof}
The explicit analytic continuation
of $I^\lam f$ is rather complicated;
see e.g. \cite{OR} and references therein
for some systematic study.

\section{Polar decomposition and Bi-Stiefel decomposition}

In this section we prove some integral
formulas related to certain polar decompositions
of matrices.
Let $S_{n, k}$ be the Stiefel manifold of
all orthonormal frames in $\mathbb K^n$. It
can be realized as the manifold
of all isometries $x\in M_{n, k}$. Let
$$G:=U(n, \bbK)=\{g\in M_{n, n}; g^\ast g=I_n\}=O(n), \, U(n),\, Sp(n)
$$
be the orthogonal, unitary, and symplectic unitary  groups according
to $\bbK =\mathbb R, \bc, \mathbb H$ respectively. 
 It
acts transitively on $S_{n, k}$ by the defining action,
and thus
$$
S_{n, k}=G/U(n-k, \bbK).
$$
Each
$x\in S_{n, k}$ defines a $k$-dimensional
subspace over $\bbK$,
$$\xi=\{x\}:=x\mathbb K^k\subset \mathbb K^n \in  G_{n, k}.
$$
 Thus
$
G_{n, k}
$ can be identified as the space of orbits in $S_{n, k}$ 
under the action of the unitary group $U(k, \bbK)$
on $\mathbb K^k$,
$$
G_{n, k}=S_{n, k}/U(k, \bbK).
$$
The manifold $G_{n, k}$ is therefore a compact Riemannian
symmetric space 
$$
G_{n, k}=U(n, \bbK)/U(n-k, \bbK)\times U(k, \bbK);
$$
see e.g. \cite{Hel1}.
 We will later
specify certain reference points in $S_{n, k}$ and  $G_{n, k}$
thus the respective isotropic subgroups. Throughout the paper
we will identify  functions $f(\xi)$  on $G_{n, k}$
with right $ U(k, \bbK)$-invariant
functions $f(x)$ on the Stiefel manifold $S_{n, k}$.

For any compact group $K$ we let $dk$ the normalized
Haar measure on $K$. Denote also 
$dv$ and $d\xi$  the normalized unique  $G$-invariant
measures on $S_{n, k}$ and on $G_{n, k}$
respectively.
The following result is  possibly known, 
for completeness we give
a proof.

\begin{lemm+}\label{polar} Let $dx$
and $dv$, $dr$ be the measures
on $M_{n, k}$, $S_{n, r}$ and $\Ome$ normalized
as above.  Almost all $x\in M_{n, k}$ can be decomposed 
uniquely as
$$
x=v r^{\frac 12}, \quad v\in S_{n,k}, \quad r\in \Ome,
$$
and under that decomposition the measure $dx$ is given by
\begin{equation*}
dx=C_0\Del(r)^{\frac d2 n} dv d\iota(r)=C_0
\Del(r)^{\frac d2 n-N/k} dv dr,
\end{equation*}
namely
$$
\int_{M_{n,k}} f(x)dx =C_0
\int_{S_{n,k}}\int_{\Ome}f(vr^{\frac 12})
\Del(r)^{\frac d2 n} dv d\iota(r),
$$
 where 
\begin{equation}\label{cons-0}
C_0=C_0(n, k)=\frac{\sqrt{\pi}^{dnk}}{\Ga_{\Ome}(\frac d{2}n)}.
\end{equation}

\end{lemm+}

\begin{proof} The polar decomposition follows from the general polar decomposition
for linear transformations and from the fact 
that the set of  elements of full rank $k$ (over
$\mathbb K$) is an open dense subset of $M_{n, k}$.
Now by the $G$-invariance we see
that there is a weight function depending only on $r$, say  $W(r)$, so that 
$$
\int_{M_{n,k}} f(x)dx  = \int_{S_{n, k}} \int_{\Ome}f(vr^{\frac 12})
W(r) d\iota(r)dv, 
$$
in term of the invariant measure $d\iota$ on $\Ome$.
To find the weight function $W(r)$ we note that for any $g\in GL_k(\bbK)$, 
  the LHS is 
$$
\int_{M_{n,k}} f(x)dx = \int_{M_{n,k}} f(xg)|\det(g)^{n}|dx
=\int_{S_{n, k}} \int_{\Ome} f(vr^{\frac 12}g)|\det(g)^{n}|W(r)
d\iota(r)dv, 
$$
where $\det(g)$ is the real Jacobian of the $\mathbb R$-linear
transformation
on $\bbK^k=\br^{dk}$, 
$y\in \br^{dk}=\bbK^k\to yg\in \bbK^k=\br^{dk}$.
We perform further  the polar decomposition of $r^{\frac 12}g$,
$$r^{\frac 12}g
=u (g^\ast r g)^{\frac 12}, \qquad u\in U(k, \bbK),
$$
so that the above integral is
$$
\int_{S_{n, k}} \int_{\Ome}
f(vu(g^\ast r g)^{\frac 12})|\det(g)^{n}|W(r)
d\iota(r)dv.
$$
Performing first  the change of variables $g^\ast r g\to r$ on $\Ome$
and then $vu\to v$ we see that it is
$$
\int_{S_{n, k}} \int_{\Ome}
f(v r^{\frac 12})|\det(g)^{n}|W((g^\ast)^{-1}rg^{-1})
d\iota(r)dv,
$$
by the invariance of $d\iota$ and respectively $dv$.
Thus the weight factor $W(r)$ transforms as 
$$
W(g^\ast r g)=|\det(g)^{n}|W(r),
$$
from which it follows, using the fact that
any $r\in\Ome$ can be written as
$r=g^\ast g$, $g\in GL_k(\mathbb K)$
(see e.g. \cite{FK-book}),  that $W(r)=c\Del(r)^{\frac{d}2 n}$ for some
constant $c$. The constant can be evaluated
by  taking the function $f(x)=e^{-\tr(x^\ast x)}$
and using the formula
    (\ref{eq:gamm}).
\end{proof}

\begin{rema+} Any $x\in M_{n,k}$, $n\ge k$, can be written as
$x=v r^{\frac 12}$ where
$r=x^\ast x \ge 0$. The element $v$ is generally not unique,
and we can always choose $v$ so that $v\in S_{n, k}$.
We mention also that,  polar decompositions, generally speaking, 
are closely related to the Bessel functions, in particular
Lemma 3.1 above is much related to  the result in \cite{FT-jfa}, p. 130.
\end{rema+}

\begin{prop+} \label{bi-stiefel} Suppose $1\le k\le k^\prime < n$
and $k+k^\prime \le n$.
\begin{itemize}
\item[(i)]
 Almost all $x\in M_{n, k}
$ can be uniquely decomposed 
as
$$
x=\begin{bmatrix} u r^{\frac 12}\\
v (I-r)^{\frac 12}
\end{bmatrix} s^{\frac 12}, \quad u\in S_{k^\prime, k},
\quad v\in S_{n-k^\prime, k}, \quad r\in (0, I), \quad s\in \Ome
$$
and under that decomposition the measure $dx$ is given by
$$
dx=C_1 \Delta(r)^{\frac d2 k^\prime - N/k}
\Delta(I-r)^{\frac d2 (n-k^\prime) - N/k} 
\Delta(s)^{\frac d2 n - N/k}dudv dr ds 
$$
with
$$
C_1=C_1(n, k^\prime, k)=\frac{\sqrt{\pi}^{d nk}}{B(\frac d2 k^\prime, 
\frac d2 (n-k^\prime))  \Ga_{\Ome}(\frac{d}{2} n)}
$$
\item[(ii)]
 Almost all $w\in S_{n, k}
$ can be decomposed 
as
$$
w=\begin{bmatrix} u r^{\frac 12}\\
v (I-r)^{\frac 12}
\end{bmatrix}, \quad u\in S_{k^\prime, k},
\quad v\in S_{n-k^\prime, k}, \quad r\in (0, I),
$$
%%s=t^2
and under that decomposition the measure $dw$ is given by
$$
dw=C_2 \Delta(r)^{\frac d2 k^\prime - N/k}
\Delta(I-r)^{\frac d2 (n-k^\prime) - N/k} 
dudv dr 
$$
with
$$
C_2=C_2(n, k^\prime)=\frac{1}{B(\frac d2 k^\prime, 
\frac d2 (n-k^\prime))}.
$$
\end{itemize}
\end{prop+}
For the proof we need another elementary result.

\begin{lemm+}\label{var-cha} The mapping
$$
(r, s)\mapsto (x, y)=(s^{\frac 12} r s^{\frac 12}, 
s^{\frac 12} (I-r)s^{\frac 12})
$$
is a diffeomorphism from $(0, I)\times \Ome$
onto  $\Ome \times \Ome $, and its Jacobian is given by
$\Del(s)^{N/k}$, namely
$$
dxdy= \Del(s)^{N/k} ds dr.
$$
\end{lemm+}
\begin{proof} We write $z=x+y=s^{\frac 12}r
s^{\frac 12} +
s^{\frac 12}(I-r)s^{\frac 12} =s$. Thus the mapping
$$(r, s )\mapsto (x, z)=(s^{\frac 12} r s^{\frac 12}, s)
$$ maps 
$\Ome\times (0, I)$ into $\{(x, z)\in \Ome\times \Ome; z >x\}$,
and its inverse is given by 
$$(x, z)\to (r, s )=(z^{-\frac 12}xz^{-\frac 12}, x).
$$ 
This proves the first part of the statement.
For the Jacobian we have  that 
$$
d\iota(z)d\iota(x)=d\iota(s)d\iota(r)
$$
by the invariance. Thus, in term of the Lebesgue measure,
\begin{equation*}
\begin{split}
dxdy&= dz dx=\Delta(x)^{N/k} \Delta(z)^{N/k} 
d\iota(x)d\iota(z)\\
&=\Delta(s)^{N/k} \Delta(r)^{N/k} d\iota(x)d\iota(z)=\Delta(s)^{N/k} dr ds.
\end{split}
\end{equation*}
\end{proof}

We prove now the Proposition.

\begin{proof}The proof is done by a change of variables.
Conceptually it is clearer to compute the integral $\int_{M_{n,k}}f(x)dx$
as in the proof of Lemma 3.1, whereas computationally it is easier
to compute the measure $dx$, and we will adopt the latter.
First of all we write $x$ as a $2\times 1$-block matrix 
under the decomposition of $\mathbb K^n=\mathbb K^{k^\prime}\oplus
\mathbb K^{n-k^\prime}$ and  then
perform polar decomposition
\begin{equation}
\label{polar-ab}
x=\begin{bmatrix} x_1 \\x_2\end{bmatrix}
=\begin{bmatrix} u_1 p_1^{1/2} \\
u_2 p_2^{1/2}  
\end{bmatrix}
\end{equation}
and we may assume that $p_1>0, p_2>0$
and $u_1\in S_{k^\prime, k}$ and $u_2\in S_{n-k^\prime, k}$, 
 by the condition 
that $k^\prime\ge k$, $n-k^\prime \ge k$. Thus 
by Lemma 3.1,
\begin{equation}
  \label{eq:dx12}
dx =dx_1dx_2,\quad 
 dx_1=c_1 du_1 \Del(p_1)^{\frac d2 k^\prime -\frac N{k}} dp_1,
\quad  dx_2=c_2 du_2 \Del(p_2)^{\frac d2 (n-k^\prime) -\frac N{k}} dp_2
\end{equation}
with $c_1=C_0(k^\prime, k)$, $c_2=C_0(n-k^\prime, k)$.
On the other hand, applying
Lemma 3.1 again to $x$ we can write
$x$ as  
$$
x=ws^{\frac 12}, \qquad w\in S_{n, k},\quad s\in \Ome
$$ 
with
\begin{equation}
  \label{eq:dw-ds}
dx=c_3 \Del(s)^{\frac d2 n -N/k}dwds, \quad c_3=C_0(n, k).
\end{equation}
We write $w\in S_{n, k}$ as a block matrix 
$w=\begin{bmatrix} w_1 \\w_2\end{bmatrix}$
and perform 
again the polar decomposition, using the fact that $w^\ast w=I=I_k$,
\begin{equation}
\label{eq:w-u-v-r}  
w=\begin{bmatrix}u r^{\frac 12}\\
v (I-r)^{\frac 12}
\end{bmatrix}
\end{equation}
for some $r$, $0\le r\le I$
and partial isometries $u$ and $v$.
We can again assume,
up to a set of $dw$-measure zero, that $0<r<I$,
$u\in S_{k^\prime, k}$ and $v\in S_{n-k^\prime, k}$.
Therefore $x$ has
the form
\begin{equation}
\label{eq:x12-rev}  
\begin{bmatrix} u_1 p_1^{1/2} \\
u_2 p_2^{1/2}\end{bmatrix} = \begin{bmatrix} x_1 \\x_2\end{bmatrix}=x=ws^{\frac 12}=\begin{bmatrix}u r^{\frac 12}\\
v (I-r)^{\frac 12}\end{bmatrix}s^{\frac 12}
=\begin{bmatrix}u r^{\frac 12}s^{\frac 12}\\
v (I-r)^{\frac 12}s^{\frac 12}\end{bmatrix}
\end{equation}
Thus
\begin{equation}
\label{eq:x12-rs}  
x_1=ur^{\frac 12}s^{\frac 12}, \quad x_2=v(I-r)^{\frac 12}s^{\frac 12};
\quad p_1=s^{\frac 12}rs^{\frac 12}, \qquad p_2=s^{\frac 12}(I-r)s^{\frac 12}.
\end{equation}

Now  using Lemma \ref{var-cha} for the change of variables
$(r, s)\to (p_1, p_2)$, we have
that  (\ref{eq:dx12}) becomes
\begin{equation*}
\begin{split}
dx&=dx_1 dx_2=c_1 c_2 
\Del(p_1)^{\frac d2 k^\prime -\frac N{k}}  \Del(p_2)^{\frac d2 (n-k^\prime) -\frac N{k}} 
du_1du_2
  dp_1  dp_2\\
&=c_1 c_2   \Del(p_1)^{\frac d2 k^\prime -\frac N{k}}  \Del(p_2)^{\frac d2 (n-k^\prime) -\frac N{k}} 
\Del(s)^{N/k}du_1du_2 ds dr\\
&=c_1 c_2   \Del(r)^{\frac d2 k^\prime -\frac N{k}} 
 \Del(I-r)^{\frac d2 (n-k^\prime) -\frac N{k}} 
\Del(s)^{\frac{d}{2}n-\frac N{K} } du_1du_2 dr ds.
\end{split}
\end{equation*}
To find the relation between $u, v$ and $u_1, u_2$
we perform further  the polar decompositions of
$r^{\frac 12}s^{\frac 12}$ and 
$(I-r)^{\frac 12}s^{\frac 12}$. They are, in view of (\ref{eq:x12-rs} ),
\begin{equation}
\label{eq:x-w-uv}  
r^{\frac 12}s^{\frac 12}=u_0 p_1^{\frac 12},
\quad (I-r)^{\frac 12}s^{\frac 12}=v_0 p_2^{\frac 12},
\end{equation}
for some $u_0, v_0\in U(k, \mathbb K)$. Thus
\begin{equation*}
u_1 p_1^{\frac 12}=x_1=ur^{\frac 12}s^{\frac 12}=uu_0p_1^{\frac 12},
 \quad u_1 p_1^{\frac 12}=x_2=v(I-r)^{\frac 12}s^{\frac 12}
=vv_0p_2^{\frac 12},
\end{equation*}
from which it follows that $u_1=uu_0, \quad u_2=v v_0$
and
$$
du_1du_2=du dv,
$$
by invariance. We now obtain
$$dx=c_1 c_2   \Del(r)^{\frac d2 k^\prime -\frac N{k}} 
 \Del(I-r)^{\frac d2 (n-k^\prime) -\frac N{k}} 
\Del(s)^{\frac{d}{2}n-\frac N{K} }du dv dr ds.
$$
This is our claim in  (i), and (ii) follows in turn by using
(\ref{eq:dw-ds}).
\end{proof}

In the next Proposition
we will  need another form of polar decomposition
that is different from the one in Lemma \ref{polar}. If $1\le k\le k^\prime$
we have that for almost all $x\in M_{k, k}$,
viewed as linear transformation
$\bbK^k \to \bbK^{k^\prime}\to 
\bbK^k$ under the natural embedding and projection,
that
$$
x=[0_{k, k^\prime-k}\,\, s^{\frac12}]u, \qquad s\in \Ome \qquad u\in S_{k^\prime, k},
$$
by the usual polar decomposition of $x$,
where $0_{p, q}$ stands for the zero $p\times q$-matrix. However
this factorization is not unique and thus no integral
formula is expected  as that in Lemma \ref{polar}. Nevertheless
we have the following substitute.
First we shall need an integral formula, proved
in \cite[Lemma 2.4]{Grinberg-Rubin}
for real Grassmannians, which follows easily
by the invariance of the measure: For any function $f$ defined on $M_{k,k}(\bbK)$
\begin{equation}\label{switch}
\int_{S_{k^\prime, k}} 
f(v^\ast u)du = \int_{S_{k^\prime, k}} f(u^\ast v )du
=\int_{S_{k^\prime, k}} f(v^\ast u )dv.
\end{equation}
\begin{prop+}\label{int-stief}
Let $1\le k\le k^\prime$.
For any  measurable function $H(w)$ on $M_{k, k}(\bbK)$,
let
$$
H_1(s)=\Del(s)^{\frac d2 (k^\prime -(k-1))-1}
\int_{S_{k^\prime, k}} H([0_{k, k^\prime-k} \,\, s^{\frac12}]v)dv, \qquad s\in \Ome.
$$
and
$$
H_2(s)=\Del(s)^{\frac d2 -1}
\int_{U(k, \bbK)} H(us^{\frac12})du, \qquad s\in \Ome.
$$
Then for any $\varepsi>-1$,
\begin{equation}\label{I-h12}
(I^{\varepsi+N/k}H_1)(s)=C_3 (I^{\frac d2(k^\prime -k)+ \varepsi+N/k}H_2)(s),
\end{equation}
where
$$
C_3(n, k^\prime, k)=\frac{C_0(k^\prime-k, k)  C_0(k,k) \Gamma_\Ome (\frac d2 (k^\prime-k)}
{C_0(k^\prime, k)}
$$
if
$k^\prime-k\ge k$
and
$$C_3(n, k^\prime, k)=\frac{C_0(k, k^\prime-k)  C_0(k,k) \Gamma_{\Ome_{k^\prime -k}} (\frac d2 k)}
{C_0(k^\prime, k)}
$$
if $k^\prime -k <k$. Here $\Gamma_{\Ome_{k^\prime -k}}$ is
the Gamma function associated with the symmetric cone $\Ome_{k^\prime -k}$
in $M_{k^\prime-k, k^\prime-k}$ and $C_0(k, k^\prime-k)$ the constant
$C_0$ in (\ref{cons-0}) with $(n, k)$ replaced by $(k, k^\prime-k)$.
\end{prop+}

\begin{proof} 
We note that
it is understood that  (\ref{I-h12}) holds whenever one
of the integrals is absolutely convergent.
The idea of the proof is to write the integral
$I^{\varepsi}H_1$ as an integral on the matrix spaces.
Fix $r\in \Ome$ and consider the integral 
 $$
\mathcal J(r):=\int_{w\in M_{k^\prime, k}:\, w^\ast w<r} H(w^\ast\begin{bmatrix} I_k\\0
\end{bmatrix}) \Del(r-w^\ast w)^{\varepsi}
dw
$$
We write the matrix $w^\ast$ also in block form:
$w^\ast=\begin{bmatrix} w_1^\ast & w_2^\ast\end{bmatrix}$ with $w_1\in 
M_{k, k}$, $w_2\in M_{k^\prime-k, k}$. We have 
$$
\mathcal J(r)=\int_{w_1 \in M_{k,k}: w_1^\ast w_1<r} H(w_1^\ast)
\left(\int_{w_2\in M_{k^\prime-k, k}: w_2^\ast w_2<r-w_1^\ast w_1}
\Del(r-w_1^\ast w_1-w_2^\ast w_2)^{\varepsi}dw_2\right)
dw_1.
$$
The inner integral can be computed 
by
$$
\int_{w_2\in M_{k^\prime-k, k}: w_2^\ast w_2<r-w_1^\ast w_1}
\Del(r-w_1^\ast w_1-w_2^\ast w_2)^{\varepsi}dw_2=c_{\varepsi} \Del(r-w_1^\ast w_1)^{\varepsi + \frac d2(k^\prime -k)}
$$
by the changing of variables $w_2=   y(r-w_1^\ast w_1)^{\frac 12}$, where
$$
c_{\varepsi}
=\int_{y\in M_{k^\prime-k, k}: y^\ast y < I_k}
\Del(I-y^\ast y)^{\varepsi} dy
$$
which will be evaluated in the end of the proof. The integral $\mathcal J(r)$
is then, by Lemma 3.1,
\begin{equation}\label{Jr-1}
\begin{split}
\mathcal J(r) &=c_{\varepsi}\int_{w_1\in M_{k, k}, w_1^\ast w_1 <r}H(w_1^\ast)
\Del(r-w_1^\ast w_1)^{\frac d2(k^\prime -k) +\varepsi} dw_1\\
&=c_{\varepsi}C_0(k, k)\int_{0}^r 
\Del(r-s)^{\frac d2(k^\prime -k) +\varepsi} 
\Del(s)^{\frac d2-1} \left( \int_{u\in U(k, \bbK)}H(us^\frac 12 ) du
\right) ds\\
&=C_{\varepsi }(I^{\frac d2 (k^\prime -k) +\varepsi+N/k}H_2)(r)
\end{split}
\end{equation}
with
$$
C_{\varepsi }= c_{\varepsi}C_0(k, k)\Gamma_\Ome(\frac d2 (k^\prime -k) +\varepsi+N/k).
$$

On the other hand, the integral $\mathcal J(r)$ can be computed
by Lemma  3.1, performing the polar decomposition $w=v s^{\frac 12}$,
 $$
\mathcal J(r)=C_0(k^\prime, k)
\int_{0}^r \Del(s)^{\frac d2(k^\prime-(k-1)) -1}
\Del(r-s)^\varepsi
\left( \int_{S_{k^\prime, k}}
H(s^{\frac 12} v^\ast \begin{bmatrix}I_k \\0\end{bmatrix})dv
\right)ds
$$
and the inner integral, by (\ref{switch}), is
$$
 \int_{S_{k^\prime, k}} H(s^{\frac 12} \begin{bmatrix}0& I_k
\end{bmatrix}v)dv, 
$$
Namely
\begin{equation*}\label{Jr-2}
\mathcal J(r)=C_0(k^\prime, k) \Gamma_{\Ome}(\varepsi+N/k)
(I^{\varepsi +N/k}H_1)(r),
\end{equation*}
with $H_1$ given as in the statement. 
Comparing the two equalities (\ref{Jr-1}) and  (\ref{Jr-2})
 we get
$$
(I^{\varepsi+N/k}H_1)(s)=C_3 (I^{\frac d2(k^\prime -k)+ \varepsi+N/k}H_2)(s),
$$
with 
$$
C_3=\frac{C_\varepsi}{C_0(k^\prime, k) \Gamma_{\Ome}(\varepsi+N/k)}.
$$

Now we evaluate $c_{\varepsi}$ and prove that
$C_3$ is as given in the Proposition (and is in fact
independent of $\varepsi$!). If $k^\prime-k\ge k$ then
by Lemma 3.1 and (\ref{eq:beta}), $c_{\varepsi} $ is
$$C_0(k^\prime-k, k)
\int_{0}^I
\Del(I-r)^{\varepsi} \Del(r)^{\frac{d}2(k^\prime-k)-N/k}dr
=C_0(k^\prime-k, k)B_\Ome(\varepsi+N/k,\frac{d}2(k^\prime-k)).
$$
If $k^\prime -k <k$, we have again,
 noticing
that for any $y\in M_{k^\prime-k, k}$, $y^\ast y < I_k$
is equivalent to $yy^\ast < I_{k^\prime -k}$, 
and that $\Del(I_k-y^\ast y)=\Del(I_{k^\prime -k}-y y^\ast)$, that
the integral can be expressed
as integration on the interval $(0, I_{k^\prime -k})$
in the symmetric cone $\Ome_{{k^\prime -k}}$ of smaller
rank,
\begin{equation*}
\begin{split}
c_{\varepsi} 
&=C_0(k, k^\prime-k) 
\int_{0<r<I_{k^\prime-k}} \Del(I_{k^\prime-k}-r)^{\varepsi} 
\Del(r)^{\frac d2 k -1-\frac d2 (k^\prime-k-1)}
dr\\
&=C_0(k, k^\prime-k) 
B_{\Ome_{k^\prime -k}}(\varepsi+1+\frac d2 (k^\prime-k-1),
\frac d2 k ).
\end{split}
\end{equation*}
The constant $C_3$ can then be computed by the formula
for $C_0$ and by the formula for the Gamma function,
we leave the elementary and yet  intricate
computations to the interested reader.
\end{proof}

%Make no difference between $\mathcal R$ and $\mathfrak R$ and

\section{Radon transform and the inverse transform}

In this section we will prove our main result,  finding
an inversion formula for the Radon transform.

To simplify notations we will write (with some abuse of
notation) $U(l)=U(l, \bbK)$ dropping
the symbol $\bbK$.
Let 
$$
x_0=\begin{bmatrix} I_k\\0 \end{bmatrix}
, \quad \hat x_0=\begin{bmatrix}0\\ I_k\end{bmatrix}\in S_{n, k}, \qquad
\xi_0=\{x_0\}, \quad \hat \xi_0=\{\hat x_0\}\in G_{n, k},
$$
be two reference points in $S_{n, k}$ and in $G_{n, k}$
respectively.
Similarly, we fix
$$
y_0=\begin{bmatrix} I_{k^\prime}\\0\end{bmatrix}
, \quad \hat y_0=\begin{bmatrix}0\\ I_{k^\prime}\end{bmatrix}\in S_{n, k^\prime}, \qquad
\eta_0=\{y_0\}, \quad \hat \eta_0=\{\hat x_0\}
\in G_{n, k^\prime}
$$
 two reference points in $S_{n, k^\prime}$ and in $G_{n, k^\prime}$.
Identifying $G_{n, k}$ with 
$G_{n, k}=G\cdot \hat \xi_0$ we have
$G_{n, k}=G/K$ where 
$$
K=U(n-k)\times U(k)
=\{k=\diag(\alpha, \delta)\in G; (\alpha, \delta)\in U(n-k)\times U(k)\}
 $$
is the isotropic subgroup of $\hat \xi_0$. Correspondingly
$$S_{n, k}=G\cdot \hat x_0 = G/U(n-k).$$
Similarly, $G_{n, k^\prime}=G\cdot \eta_0=G/K^\prime$,
$S_{n, k^\prime}=G\cdot y_0=G/U(n-k^\prime)$,
with
$$
K^\prime=U(k^\prime)\times U(n-k^\prime)
=\{k=\diag(\tau, \rho)\in G; (\tau, \rho)\in U(k^\prime)\times U(n-k^\prime)\} $$
the isotropic subgroup of $\eta_0$.
We will henceforth
 fix these realizations of the groups
$U(k)$, $U(n-k)$, $U(k^\prime)$
and $U(n-k^\prime)$ when viewed as subgroups of $G=U(n)$,
if no ambiguity would arise.

For any  $\eta\in G_{n, k^\prime}$ the subset
$$
S(\eta)=\{\xi\in G_{n, k}; \xi\subset \eta\}
$$
is a totally geodesic submanifold of $G_{n, k}$. It is itself
a symmetric space with the induced metric. For $\eta=\eta_0$,
$$
S(\eta_0)=G_{k^\prime, k}=U(k^\prime)/U(k)\times U(k^\prime -k)
$$
where $U(k)\times U(k^\prime -k)$
(now with a different realization) consists of the elements in $G$
of the form $\diag(\delta, \epsilon, I_{n-k^\prime})$. For any
$\eta\in G_{n, k^\prime}$  let $\eta=g_\eta \eta_0$ with $g_\eta\in G$, then
$$
S(\eta)=g_\eta S(\eta_0).
$$
We let $d_{\eta}\xi$ be the unique $U(k^\prime)$ invariant
measure on $S(\eta)$ via the above identification.

We define the Radon transform
$\mathcal R: C^\infty(G_{n, k})\to C^\infty(G_{n, k^\prime})$ by
\begin{equation}\label{def-radon-2}
(\mathcal Rf)(\eta)=\int_{S(\eta)} f(\xi)d_\eta\xi.
\end{equation}
Equivalently it can be defined as
\begin{equation}\label{def-radon}
(\mathcal Rf)(\eta)=\int_{U(k^\prime)} f(g_\eta \tau \xi_0)d\tau.
\end{equation}
Similarly we define the Radon transform on 
$S_{n,k}$, 
\begin{equation}\label{def-radon-var}
(\mathcal Rf)(y)=\int_{U(k^\prime)} f(g_y \tau x_0)d\tau.
\end{equation}
It is easy to see that $\mathcal Rf$ is well-defined
and maps $C^\infty(S_{n, k})$
to $C^\infty(S_{n, k^\prime})$. In particular the two
formulas agree when acting on $f\in C^\infty (G_{n, k})$
viewed as right $U(k)$-invariant functions 
on $S_{n, k}$.

As explained in the introduction, we shall define  the operator 
$\phi(\eta)\to (\mathcal T_{r^{1/2}} \phi)(\xi) $ 
from $C^\infty(G_{n, k^\prime})$ to $C^\infty(G_{n, k})$, 
as the integration of functions $\phi$
over planes $\eta$ that are of angle $(x, \eta)$ so
that $\Cos^2(x, \eta)=u ru^\ast$, $u\in U(k)$, with
$x\in S_{n, k}$ a representative of $\xi$; we shall
define a slight generalization, $T_a$, for a $a\in M_{k, k}$, $a^\ast a\le I$.
 We will write the 
integration as one  on  the group $U(n-k)$.
For that purpose 
we define the following  $n\times n$ matrices, 
written in block matrices
under the decomposition $\mathbb K^n=\bbK^{k^{\prime}-k}
\oplus\bbK^{k}\oplus
\bbK^{n-k^{\prime}-k}
\oplus\bbK^{k}$,
$$
j(a)=\begin{bmatrix}I_{k^\prime-k} &0 &0 &0\\
0& (I_{k}-a^\ast a)^{\frac12} &0 &a^\ast\\
0& 0 & I_{n-k^\prime -k} &0 \\
0& -a &  0 &(I_{k}-a a^\ast)^{\frac12}
\end{bmatrix},
$$
and
$$
h(a)=\begin{bmatrix}
a^\ast &0 &0 & (I_{k}-a^\ast a)^{\frac12}\\
0& I_{n-k^\prime -k} &0 &0\\
0& 0& I_{k^\prime-k} & 0 \\
-(I_{k}-a a^\ast)^{\frac12} & 0 & 0& a
\end{bmatrix}.
$$
where $a\in M_{k, k}$, $a^\ast a\le I$.

Some simple observations
 about those elements
are given in the next Lemma.
\begin{lemm+} 
\begin{itemize} 
\item[(i)]
Suppose $0\le a^\ast a \le I_k$.
The elements $j(a)$ and $h(a)$ are in $G=U(n, \mathbb K)$,
and $j(a)^{-1}=j(a)^\ast =j({-a})$,
$h(a)^{-1}=h(a)^\ast =h({a^\ast})$. 
\item[(ii)]The  projection $P_{j(a)\eta_0} =j(a)P_{\eta_0}j(a)^\ast$
 onto the subspace $j(a)\eta_0=\{j(a)y_0\}\in G_{n, k^\prime}$ is
given by
$$
P_{j(a)\eta_0} = \begin{bmatrix}
I_{k^\prime-k} &0   &0    & 0\\
0              &I_k -a^\ast a & 0 & (I_{k}-a^\ast a)^{\frac12} a^\ast\\
0& 0& 0& 0 \\
 0 &a(I_{k}-a^\ast a)^{\frac12} & 0& aa^\ast
\end{bmatrix}.
$$

\item[(iii)]
For any $v\in S_{k^\prime, k}$, we have
$$
j(a)^{-1}\begin{bmatrix} v
\\ 0\end{bmatrix}= \begin{bmatrix} \begin{bmatrix}I_{k^\prime-k} &0\\
0& (I_{k}-a^\ast a)^{\frac12} \end{bmatrix} v\\
0_{n-k^\prime-k, k}\\
\begin{bmatrix} 0&a\end{bmatrix} v
\end{bmatrix},
$$
and
$$
h(a)\begin{bmatrix}  0
\\ v\end{bmatrix}
=\begin{bmatrix} \begin{bmatrix} 0_{k, k} 
&(I_{k}-a^\ast a)^{\frac 12}\end{bmatrix}v\\
0_{n-k^\prime-k}\\
\begin{bmatrix} I_{k^\prime-k} &0\\
0& a\end{bmatrix} v
\end{bmatrix}.
$$
\item[(iv)] Assume that $0<aa^\ast <I_k$. For any $v\in S_{k^\prime, k}$ there exists
an $l\in U(n-k)$ such that
\begin{equation}
\label{ia-v-hb}
j(a)^{-1}\begin{bmatrix} v\\
0_{n-k^\prime, k}\end{bmatrix}
=\diag(l, I_{k})h({b})\hat x_0
\end{equation}
where $b=\begin{bmatrix} 0_{k, k^\prime-k}&a\end{bmatrix} v\in M_{k, k}$. 
\end{itemize} 
\end{lemm+}

\begin{proof} We prove only the last statement (iv), the remaining
 are proved by simple matrix computations.
First we see that $h(b)\hat x_0$
is the right hand side of (\ref{ia-v-hb}) is,
by (iii),
\begin{equation}\label{h-on-hat x} 
h(b)\hat x_0 = \begin{bmatrix} 
(I_{k}-a^\ast a)^{\frac 12}\\
0_{n-2k, k}\\
a 
\end{bmatrix}.
\end{equation}
We have  the left hand side of (\ref{ia-v-hb}),
by the first formula in (iii), is of the form
$$
j(a^{-1})\begin{bmatrix} v\\
0_{n-k^\prime, k}\end{bmatrix}
=\begin{bmatrix} q \\b\end{bmatrix}\in S_{n, k}
$$
for some $q\in M_{n-k, k}$, $q^\ast q= I_k - b^\ast b$.
By our assumption $0<a^\ast a<I$ we see
that  $b=[0 \, \,a]v$ satisfies also $b^\ast b <I$
and which in turn implies that $q$ is of full rank $k$.
The polar decomposition
of $q$ is then 
$q=u (I_k -b^\ast b)^{\frac 12}$,
with $u\in S_{n-k, k}$, and $u$ can
be further written as
$$
u=l  \begin{bmatrix} I_k \\0_{n-k-k, k}\end{bmatrix}
$$ 
for
$l\in U(n-k)$.
 Namely
$$
j(a^{-1})\begin{bmatrix} v\\
0_{n-k^\prime, k}\end{bmatrix}
=\diag(l, I_{k})\begin{bmatrix}(I_{k}-b^\ast b)^{\frac 12}\\0
  \\b\end{bmatrix}\in S_{n, k},
$$ 
which is the right hand side
(\ref{ia-v-hb}).
\end{proof}

We define, for any $b\in M_{k,k}(\mathbb K)$, $b^\ast b\le I$, 
the operator $\mathcal T_{b}: C^\infty(S_{n, k^\prime})
\to C^\infty(S_{n, k})$ 
\begin{equation}\label{def-T-a}
\mathcal T_{b}
\phi(x)=\int_{U(n-k)\times U(k)
}\phi(g_x k j(b)^{-1} y_0)dk, \quad x=g_x \hat x_0, \quad g_x\hat\in G.
\end{equation}
It maps $C^\infty(G_{n, k^\prime})
$ to $C^\infty(G_{n, k})$ and has the form
\begin{equation}\label{def-T-a-var}
\mathcal T_{b}
\phi(\xi)=\int_{U(n-k)\times U(k)
}\phi(g_\xi k j(b)^{-1} \eta_0)dk.
\end{equation}

The following lemma clarifies
the  geometric
meaning of the above integral.
For completeness we give an elementary
proof.
\begin{lemm+}
Let $0\le r \le I=I_k$ and $\eta\in G_{n, k^\prime}$ and
 $x=g_x\hat x_0\in S_{n, k}$ for $g_x\in G$.
Then $\Cos^2(\eta, x)=\delta_0 r\delta_0^\ast$ for
some $\delta_0\in U(k)$ if and only
if $\eta=g_x\rho j({r^\frac 12}) \eta_0$ for some
$\rho =\diag (\alpha, \delta)\in K$.
\end{lemm+}
\begin{proof} Suppose $\rho =\diag (\alpha, \delta)\in K$ and
 $\eta=g_x\rho j({r^\frac 12})^{-1} \eta_0$.
The projection $P_{\eta}\in M_{n, n}$
is given by
$$
P_{\eta}g_x\rho  P_{j({r^\frac 12})^{-1} \eta_0}
\rho^\ast g_x^\ast 
$$
and
$$
\Cos^2(\eta, x)=x^\ast P_{\eta} x
=\hat {x_0}^\ast
\rho P_{j({r^\frac 12})^{-1} \eta_0}
\rho^\ast \hat {x_0}= \delta r \delta^\ast
$$
by Lemma 4.1 (ii).

Conversely, suppose $\Cos^2(\eta, x)=\delta_0 r \delta^\ast:=s$.
We let $\ome=g_x^{-1}\eta$. We prove that $\ome=\rho j(r^{\frac 12})^{-1} \eta_0$
for some $\rho\in K$,
which is our claim. Now $P_{\eta}=g_x P_{\ome} g_x^\ast$,
and then  by definition
$$s=\Cos^2(\eta, x)=x^\ast P_{\eta}x = 
\hat x_0^\ast P_{\ome} \hat x_0.
$$
 The matrix $P_{\omega}$ is a self-adjoint matrix on
$\bbK^n$, and under the decomposition of $\bbK^n=\bbK^{n-k}\oplus
\bbK^{k}$ it is of the block-matrix form
$$
P_{\omega}=\begin{bmatrix} A& B\\
B^\ast & s\end{bmatrix}.
$$
by the previous formula,
with $A$ self-adjoint and $B\in M_{n-k, k}(\bbK)$.
Being a projection, $P^2_{\ome}=P_\ome$, that is
$$
A^2 +BB^\ast =A, \quad AB+Bs=B, \quad B^\ast B+ s^2= s.
$$
The  polar decomposition of $B$ is of the form,
under the decomposition $\mathbb K^{n-k}=\mathbb K^{k^\prime-k}\oplus\mathbb K^{k}
\oplus\mathbb K^{n-k-k^\prime}
$,
$$
B=\alpha_1 \begin{bmatrix} 0\\
s^{\frac 12}(I_k -s)^{\frac 12}\\
0\end{bmatrix}, \quad \alpha_1\in U(n-k),
$$
for $B^\ast B=s(I_k-s)$.

Let $A_1=\alpha_1^{-1}A\alpha_1$. Writing $A_1$ in block
matrix form under the same decomposition
and using the equation $A^2+BB^\ast =A$ we see that
$$
A_1= \begin{bmatrix} A_{11}& 0 &A_{13}\\
0&(s)^{\frac 12}(I_k -s)^{\frac 12}&0
\\A_{13}^\ast &0&A_{33}\end{bmatrix}.
$$
and that the matrix
$$
\begin{bmatrix} A_{11}& A_{13}\\
A_{13}^\ast &A_{33}\end{bmatrix}.
$$
is a projection of rank $k^\prime -k$ since $P_\ome$ is of rank $k^\prime$.
Thus 
$$
\begin{bmatrix} A_{11}& A_{13}\\
A_{13}^\ast &A_{33}\end{bmatrix}
=\alpha_2\begin{bmatrix} I_{k^\prime-k}& 0\\
0 &0\end{bmatrix}\alpha_2^\ast
$$
for some $\alpha_2\in U(n-2k)$.
Finally we
have
$$
P_\omega =\rho j(r^{\frac 12})^{-1} P_{\eta_0} j(r^{\frac 12}) \rho^\ast
$$
with $\rho=\diag(\alpha_1 \alpha_2, \delta)\in K$,
which  is equivalent to
$$
\omega=\rho j(r^{\frac 12})^{-1} \eta_0,
$$
as any subspace $\omega$ is uniquely determined by the projection
$P_\omega$.
\end{proof}

For $0 <r < I$
the operator $\mathcal T_{r^\frac 12}$
can be
expressed in terms of the angle explained in the introduction,
$$
(\mathcal T_{r^\frac 12}\phi)(x) =\int_{\eta=\{y\}: \, \Cos^2(\eta, x)=uru^\ast, 
u\in U(k)}\phi(y) d_x y.
$$
For the real Grassmannians this is the slightly
corrected version of the formulas (1.8) and (3.2) in \cite{Grinberg-Rubin}
($uru^\ast$ should be in place of $r$ there).

To state the next proposition we need also a mean-value operator,
$\mathcal W_{b}: C^\infty(S_{n, k})\to 
C^\infty(S_{n, k})$, for $b\in M_{k, k}$, $b^\ast b\le I$,
\begin{equation}\label{operator-w}
(\mathcal W_{b}f)(x)=\int_{U(n-k)}f(g_x\diag(\alpha, I_k) h(b) \hat x_0)d\alpha;
\end{equation}
again $\mathcal W_a$  is well-defined. Observe that it does
not map functions on $G_{n, k}$ to functions $G_{n, k}$.
If $f\in C^\infty(G_{n, k})$,
the function
$$
\int_{U(k)}(\mathcal W_{\delta a}f)(x)d\delta
$$
is right $U(k)$-invariant and thus defines
a functions on $G_{n, k}$. 
(For $k=1$ the operator
$\mathcal W_{\delta a}$ can roughly  speaking
be incorporated into the operator $\mathcal T$
so that it doesn't play a role.)

First we state an integral formula, which
is a direct consequence of Proposition 3.3 (ii), written in the form
$$
w=\begin{bmatrix} v (I-r)^{\frac 12}
\\
u r^{\frac 12}
\end{bmatrix}, \quad v\in S_{n-k, k},
\quad u\in U(k), \quad r\in (0, I),
$$
and of the
formula (\ref{h-on-hat x});
see also \cite[Lemma 3.5]{Grinberg-Rubin}
(a different formula for the mean-value
operator is defined there).

\begin{lemm+}\label{int-S-nk-W}
Suppose $1\le p\le \infty$ and $f\in L^p(S_{n, k})$. Then
for any $x\in S_{n,k}$
\begin{equation*}
\begin{split}
&\quad\, \int_{S_{n,k}} f(g_x w)dw=\int_{S_{n,k}} f(w)dw\\
&=C_2(n, k)\int_0^{I_k} \Del(I_k - r)^{\frac{d}{2}(n-2k+1)-1}
 \Del(r)^{\frac{d}{2}-1}
\left( \int_{U(k)} \big(\mathcal W_{u r^{\frac 12}}f\big)(x) du\right) dr.
\end{split}
\end{equation*}
In particular, the function 
\begin{equation}\label{Phi0}
\Phi_0(r)= \Phi_0(r, x):=\Del(r)^{\frac{d}{2}-1}
\int_{U(k)} (\mathcal W_{u r^{\frac 12}}f\big)(x)du
\end{equation}
is in $L^1((0, I_k), \Del(I_k - r)^{\frac{d}{2}(n-2k+1)-1}
dr)$.
\end{lemm+}

\begin{prop+}Let $\varepsi\in \mathbb C$. Suppose $f\in L^1(S_{n, k})$ and $\phi=\mathcal R f$. Then, in the sense of distribution,
$$
I^{\varepsi +N/k}\Phi_1
=C_{3}I^{\frac d2 (k^\prime -k)+\varepsi +N/k} \Phi_0,
$$
as analytic continuation in $\varepsi$,
where 
$$\Phi_1(s)=\Phi_1(s, x)=\Delta(s)^{ \frac d2(k^\prime-(k-1))-1} (\mathcal T_{s^{\frac 12}}\phi)(x)
$$
and 
$\Phi_0$ is as in (\ref{Phi0}).
\end{prop+}
\begin{proof} Let $r\in (0, I)$. We compute $(\mathcal T_{r^{\frac 12}}\phi) (x)
=(\mathcal T_{r^{\frac 12}}\mathcal Rf)(x)$. We prove
the formula for $\varepsi\ge 0$,  in which
case the  convergence
of the relevant integrals are easily justified
by using Lemma \ref{int-S-nk-W}, and the result then follows for any $\varepsi$
by analytic continuation.

 By the definitions of $\mathcal T$ and $\mathcal R$ we have, 
$$
(\mathcal T_{r^{\frac 12}}\phi) (x)
=\int_{K} \int_{\tau\in U(k^\prime)} f(g_x \rho j(r^{\frac 12})^{-1} \widetilde\tau x_0 )
d\tau d\rho,
$$
where $\tilde\tau =\diag(\tau, I_{n-k})$ in the formula.
Introducing the variable $v= \tau I_k\in S_{k^\prime, k}$, we have
$$
\widetilde \tau x_0=\begin{bmatrix} v\\ 0\end{bmatrix};
$$
 furthermore using  Lemma 4.1 (iv)
we have 
$$
j(r^{1/2})^{-1}\widetilde \tau x_0=\diag(l, I_k) h([0\, \, r^{\frac 12}] v)\hat x_0,
$$
for some $l\in U(n-k)$.
Now  $\rho=\diag (\alpha, \delta)
$, 
$$
\rho j(r^{1/2})^{-1}\widetilde \tau x_0=\diag (\alpha, \delta)\diag(l, I_k) h([0\, \, r^{\frac 12}] v)\hat x_0
=\diag (\alpha l, \delta)h([0\, \, r^{\frac 12}] v)\hat x_0,
$$
thus changing variables $\alpha l \to \alpha$
the previous integral becomes
\begin{equation*}
\begin{split}
&\quad\,\int_{S_{k^\prime, k}} \int_{\alpha \in U(n-k))}\int_{\delta\in U(k)}
 f(g_x \diag (\alpha, I_k)\diag(I_{n-k}, \delta)h([0\, \, r^{\frac 12}] v)\hat x_0) )
d\delta d\alpha dv
\\
&=\int_{S_{k^\prime, k}} \int_{\alpha \in U(n-k))} \int_{\delta\in U(k)}
 f(g_x \diag (\alpha, I_k)h(\delta[0\, \, r^{\frac 12}] v)\hat x_0) )d\delta d\alpha dv.
\end{split}
\end{equation*}
In term of the mean-value operator $\mathcal W_b$ 
the integral
can be written as
\begin{equation}
(\mathcal T_{r^{\frac 12}}\phi) (x)
=\int_{S_{k^\prime, k}} \int_{\delta\in U(k)}\big(\mathcal W_{\delta [0, r^{\frac 12}]v}f\big)(x)
d\delta dv
\end{equation}
We use now Proposition 3.5 with $H(w)=\int_{\delta\in U(k)}
(\mathcal W_{\delta w}f)(x)$,
obtaining
$$
(I^{\varepsi+N/k} \Phi_1) (s)
=C_3 (I^{\frac d2(k^\prime -k)+\varepsi+N/k}\Phi_0)(s), \qquad s\in (0, I).
$$
\end{proof}

The next lemma
is elementary, the first two parts can be proved by the same method
as that of Lemma 3.4 in \cite{Grinberg-Rubin}. The last part is proved
by using the triangle inequality and the fact that $\Del(s)\le 1$ for $s\in(0, I)$;
 we omit the proof.
\begin{lemm+} \label{approx-id} 
\begin{itemize}
\item [(i)]
 For any $f\in L^p(S_{n, k})$, $1\le p<\infty$, we have
$$
\lim_{a\to I_k, a^\ast a<I_k} \mathcal W_{a}f =f.
$$
where the limit is in $L^p(S_{n, k})$-sense. 
\item [(ii)] 
If $f\in L^p(G_{n,k})$ then
$$
\lim_{a\to I_k, a^\ast a<I_k} \int_{U(k)}(\mathcal W_{\delta a}f) d\delta =f.
$$
and the limit is taken in $L^p(G_{n, k})$.
\item [(iii)]  If  $f\in L^p(S_{n,k})$ or $L^p(G_{n,k})$ then for any  complex number
$\zeta$,
$$
\lim_{a\to I_k, a^\ast a<I_k} \Del(a^\ast a)^\zeta\int_{U(k)}
(\mathcal W_{\delta a}f) d\delta =f.
$$
and the limit is taken in the respective spaces.
\end{itemize}
\end{lemm+}

\begin{theo+} Let $\mathbb K=\mathbb R, \bc, \mathbb H$
be the field of real, complex and quaternionic numbers.
Let $1\le k < k^\prime< n-1$, $k+k^\prime\le n$.
Let $\Del(\partial)$ be the Cayley-type differential
operator acting on the space of self-adjoint
$k\times k$-matrices over $\bbK$ and $I^{\lambda}$ be the
G\aa{}rding-Gindikin fractional integration.
Suppose    $1\le p<\infty$,
$f\in L^p(G_{n, k}(\mathbb K)) $ and $\phi=\mathcal R f$.
Let $m$ %$m=\frac a2 +\varepsi$
be such that
$$
m>\frac d2 (k^\prime-1).
$$
Then the Radon transform $f\to \phi=\mathcal Rf$ is inverted by
the formula
$$
f=\lim_{s\to I_k }\Del(\partial_s)^{m}
I^{m-\frac d2(k^\prime -k)} \Phi$$
where the limit is taken in the space $L^p(G_{n, k}(\bbK))$ and the
differential operator acts in the sense of distribution,  and where
$$\Phi(s)= \Phi(s, \xi)=\Delta(s)^{\frac d2(k^\prime-(k-1))-1}\big(\mathcal T_{s^{\frac 12}}
\phi\big)(\xi).
$$
 \end{theo+}
 \begin{proof} Let $m$ be as in the statement and  let $
\varepsi>-1$ be given so that $m=\frac d2(k^\prime -k) +\varepsi +N/k$.
Then Proposition 4.4 can be restated as
$$
(I^{m} \Phi_0)(s)
=(I^{m-\frac d2(k^\prime -k}\Phi )(s).
$$
Now Lemma 4.3 implies that $\Phi_0$ is in
a $L^1$-space, and we can then use
Lemma 2.1 (iii), obtaining
$$
 \Del(s )^{\frac d2 -1} \int_{U(k)}
(\mathcal W_{u s ^{\frac 12}}f)(x) d\delta
=\Phi_0(s)=\Del(\partial_s)^m (I^{m-\frac d2(k^\prime -k}\Phi )(s),
$$
in the distributional sense.
The result then follows from Lemma \ref{approx-id} (iii).
 \end{proof}

There arise several interesting questions
such as inverting Radon transform
on non-compact symmetric matrix domains and
finding a Plancherel formula.
In a forthcoming paper we will generalize
the results in this paper to that setup.

\newcommand{\noopsort}[1]{} \newcommand{\printfirst}[2]{#1}
  \newcommand{\singleletter}[1]{#1} \newcommand{\switchargs}[2]{#2#1}
  \def\cprime{$'$} \def\cprime{$'$}
\providecommand{\bysame}{\leavevmode\hbox to3em{\hrulefill}\thinspace}
\providecommand{\MR}{\relax\ifhmode\unskip\space\fi MR }
\providecommand{\MRhref}[2]{%
  \href{http://www.ams.org/mathscinet-getitem?mr=#1}{#2}
}
\providecommand{\href}[2]{#2}

\end{document}